\documentclass{aptpub}
\authornames{P.W.~GLYNN AND Z.~ZHENG} 
\shorttitle{Approximation of Markov Chain Expectations} 

\DeclareMathOperator{\E}{\mathbb{E}}

\usepackage{amsmath}
\usepackage{mathtools}

\usepackage{tikz}

\begin{document}

\title{Approximation of Markov Chain Expectations and the Key Role of Stationary Distribution Convergence} 

\authorone[Stanford University]{Peter W. Glynn} 
\authorone[University of California Berkeley]{Zeyu Zheng}


\addressone{475 Via Ortega, Stanford, CA 94305; 4125 Etcheverry Hall, Berkeley, CA 94720;} 
\emailone{glynn@stanford.edu; zyzheng@berkeley.edu} 

\begin{abstract}

Consider a sequence $P_n$ of positive recurrent transition matrices or kernels that approximate a limiting infinite state matrix or kernel $P_{\infty}$. Such approximations arise naturally when one truncates an infinite state Markov chain and replaces it with a finite state approximation. It also describes the situation in which $P_{\infty}$ is a simplified limiting approximation to $P_n$ when $n$ is large. In both settings, it is often verified that the approximation $P_n$ has the characteristic that its stationary distribution $\pi_n$ converges to the stationary distribution $\pi_{\infty}$ associated with the limit. In this paper, we show that when the state space is countably infinite, this stationary distribution convergence implies that $P_n^m$ can be approximated uniformly in $m$ by $P_{\infty}^m$ when n is large. We show that this ability to approximate the marginal distributions at all time scales $m$ fails in continuous state space, but is valid when the convergence is in total variation or when we have weak convergence and the kernels are suitably Lipschitz. When the state space is discrete (as in the truncation setting), we further show that stationary distribution convergence also implies that all the expectations that are computable via first transition analysis (e.g. mean hitting times, expected infinite horizon discounted rewards) converge to those associated with the limit $P_{\infty}$. Simply put, we show that once one has established stationary distribution convergence, one immediately can infer convergence for a huge range of other expectations.

\end{abstract}



\section{Introduction}

Let $S$ be a finite or countably infinite state space, and for $1
\le n \le \infty$, let $P_n = (P_n(x, y) : x, y \in S)$ be a one-step transition matrix indexed by $S$. Suppose that for $1\le n \le \infty$, $P_n$ possesses a unique stationary distribution $\pi_n$. If, for each $x, y \in S$,
\begin{equation}
    P_n(x, y) \to P_\infty(x, y)  \label{eq:1.1}
\end{equation}
as $n \to \infty$, this raises the natural question of when
\begin{equation}
    \pi_n(y) \to \pi_\infty(y) \label{eq:1.2}
\end{equation}
as $n \to \infty$, for each $y \in S$. When $P_n$ is irreducible and aperiodic for $1\le n \le \infty$,
\begin{equation}
    P_n^m(x, y) \to \pi_n(y) \label{eq:1.3}
\end{equation}
as \( m \to \infty \), and \eqref{eq:1.2} and \eqref{eq:1.3} assert that for each \( x, y \in S \), the iterated limits commute, so that
\begin{equation}
    \lim_{m \to \infty} \lim_{n \to \infty} P_n^m(x, y) = \lim_{n \to \infty} \lim_{m \to \infty} P_n^n(x, y),\label{eq:1.4}
\end{equation}
We refer to the interchangeability of the limits in \( m \) and \( n \) associated with \eqref{eq:1.4} as the \textit{weak interchange property}  and to \eqref{eq:1.2} as the \textit{stationary distribution convergence property}. The weak interchange property implies that whether we proceed vertically and then horizontally in the diagram below, versus horizontally and then vertically, we get the same result.
\[
\begin{tikzpicture}[node distance=2cm, auto]
    \node (Pnm) {\(P_n^m\)};
    \node (P) [right of=Pnm, node distance = 4cm] {\(P_\infty^m\)};
    \node (pin) [below of=Pnm] {\(\pi_n\)};
    \node (pi) [below of=P] {\(\pi_\infty\)};

    \draw[->] (Pnm) to node {\(n \to \infty\)} (P);
     \draw[->] (P) to node {\(m \to \infty\)} (pi);
    \draw[->] (Pnm) to node [swap] {\(m \to \infty\)} (pin);
    \draw[->] (pin) to node [swap] {\(n \to \infty\)} (pi);
    \draw[->, dashed] (Pnm) to node [swap] {} (pi);
\end{tikzpicture}
\]

But this raises the question of the dotted line \textit{diagonal convergence}. In particular, when is it true, in the presence of \eqref{eq:1.1}, \eqref{eq:1.2}, and \eqref{eq:1.3}, that whenever \( m_n \to \infty \), we have that
\begin{equation}
    P_n^{m_n}(x, y) \to \pi_\infty(y) \label{eq:1.5}
\end{equation}
as $n\to \infty$, for  $x, y \in S$?

When \eqref{eq:1.5} holds for every sequence \( m_n \to \infty \), this is equivalent, in the presence of \eqref{eq:1.1}, \eqref{eq:1.2}, and \eqref{eq:1.3}, to requiring that
\begin{equation}
    \sup_{m\ge 0} \sum_{y \in S} |P_n^m(x, y) - P^m_\infty(x, y)| \to 0. \label{eq:1.6}
\end{equation}
as $n\to\infty$. We refer to \eqref{eq:1.6} as the \textit{strong interchange property}, since it implies that the marginals of the Markov chain driven by \( P_n \) can be well approximated by those of \( P_\infty \) at all time scales, namely short, long, and intermediate scales. (The weak interchange property guarantees only that the marginals of the chain \( (X_n:{n \geq 0} )\) under \( P_n \) are well approximated for short scales (fixed \( m \)) and long scales (in equilibrium).) Note that the strong interchange property ensures, for example, that all lag-$m$ covariances between \(f(X_j)\) and \(f(X_{j+m})\) computed under \(P_n\) are uniformly (in \(j\) and \(m\)) close to those computed under \(P_\infty\), so that the mixing structure can be suitably approximated (at least when \(f\) is bounded).

Perhaps surprisingly, the weak interchange property implies diagonal convergence and the strong interchange property in this setting of countable state spaces; see Corollary 1. However, in Section 4 of this paper, we provide a counter-example when \( S = [0, 1] \) in which the convergence in \eqref{eq:1.1}, \eqref{eq:1.2}, and \eqref{eq:1.3} is replaced by weak convergence. While the weak interchange and strong interchange properties are generally not equivalent in continuous state space, we prove that they are equivalent when \eqref{eq:1.1}, \eqref{eq:1.2}, and \eqref{eq:1.3} are generalized to chains exhibiting total variation convergence (Proposition 1), and when weak convergence is supplemented by a suitable Lipschitz inclusion condition on the transition kernel of the chain (Proposition 2). When Proposition 2 is applied to the delay sequence of the \( G/G/1 \) queue, we obtain a generalization of a result of \cite{breuer2008continuity}; see Example 2.1.

The weak interchange problem has been studied in the setting of diffusion approximations for queues; see \cite{gamarnik2006validity}, \cite{budhiraja2009stationary}, \cite{gurvich2014diffusion}, \cite{ye2016diffusion}, \cite{braverman2017stein}, \cite{lee2020stationary}, \cite{ye2021diffusion} and \cite{dai2023asymptotic}. This issue also arises in the study of continuity of stochastic models as a function of their ``input sequences"; see, for example, \cite{kennedy1972continuity,whitt1974continuity, whitt1980some}. The reference \cite{karr1975weak} provides conditions guaranteeing that a sequence of continuous state space Markov chains exhibits the weak interchange property.

As noted above, conditions \eqref{eq:1.1} to \eqref{eq:1.3} ensure that the Markov chain associated with \( P_n \) is well approximated at short time scales and long time scales. When \( S \) is countably infinite, stationary distribution convergence implies diagonal convergence. One might hope that virtually all performance measures for the chain evolving under \( P_n \) can be approximated by the corresponding performance measure under \( P_\infty \). This turns out to be true in the countable state setting; see Theorem 2 and Proposition 3.

We conclude by noting that the questions discussed in this paper are also highly relevant to the extensive literature on numerical truncation methods for Markov chains. Early work includes \cite{seneta1967finite}, \cite{golub1973computation}, \cite{golub1974computation}, \cite{seneta1980computing}, and \cite{wolf1980approximation}, in which conditions are developed that ensure that the weak interchange property holds for a sequence \( P_n \) of transition matrix approximations to \( P_\infty \). This is relevant to numerical algorithms intended to compute the stationary distribution of the limiting infinite state transition matrix \( P_\infty \). Similarly, one may wish to use the finite dimensional matrix \( P_n^m \) as an approximation to \( P_\infty^m \), or to compute expectations for \( P_\infty \) involving first transition analysis by solving corresponding finite-dimensional linear systems involving \( P_n \). The results in the current paper make clear that any truncation method that preserves stationary distribution convergence is guaranteed to also consistently approximate all these other performance measures associated with \( P_\infty \). This is especially useful, in view of the fact that the great majority of the existing literature on state space truncation has focused on the question of stationary distribution convergence.

This paper is organized as follows. Section 2 develops a general theorem (Theorem 1) that provides conditions under which the weak interchange property implies the strong interchange property, and then specializes the result first to the total variation convergence setting (Proposition 1) and then to the weak convergence setting (Proposition 2). Section 3 shows that when the state space is countably infinite, the stationary distribution convergence property implies convergence of all the expectations that are computable via first transition analysis. Section 4 provides a discussion of a counter-example demonstrating that the weak interchange property does not imply the strong interchange property in general, when the Markov chain has a continuous state 
space. Finally, Section 5 extends the theory to the Markov jump process settings. 

\section{The Weak and Strong Interchange Properties}

Let \( X = (X_m:{m \geq 0} )\) be an \( S \)-valued Markov chain. (We do not assume in this section that \( S \) is finite or countably infinite, unless otherwise stated.) For \( 1 \leq n \leq \infty \) and \( x \in S \), let \( P_{n,x}(\cdot) \) (\( E_{n,x}(\cdot) \)) be the probability (expectation) on the path-space \( \Omega = S^\infty \) under which \( X_0 = x \) and
\[
P_{n,x}(X_{m+1} \in dz \mid X_m = y) = P_n(y, dz)
\]
for \( y, z \in S \) and \( m \geq 0 \), where \( P_n = ( P_n(y, dz):y,z \in S ) \) is a transition kernel on \( S \).

Given a finite signed measure \( \nu \) on \( S \), and a set \( \mathcal{F} \) of (suitably measurable) bounded real-valued functions, put
\[
\| \nu \|_{\mathcal{F}} \triangleq \sup_{f \in \mathcal{F}} |\int_S f(x)\nu(dx)|.
\]

We note that, in general, \( \| \cdot \|_{\mathcal{F}} \) is not a norm, since \( \| \nu \|_{\mathcal{F}} = 0 \) may not imply that \( \nu = 0 \). When \( \mathcal{F}_w \) is the set of  functions \( f \) for which \( |f(x)| \leq w(x) \) for \( x \in S \) (with \( w \) positive), we denote \( \| \nu \|_{\mathcal{F}_w} \) as \( \| \nu \|_w \).

\textbf{Remark 2.1.} The quantity \( \|\cdot\|_w \) is the \( w \)-weighted total variation norm. When \( w = e \) with \( e(x) = 1 \) for \( x \in S \), this corresponds to the (standard) \textit{total variation} norm. Given a sequence \( (\eta_n:{1\le n \le \infty}) \) of probabilities for which \( \|\eta_n - \eta_\infty\|_e \to 0 \) as \( n \to \infty \), we write \( \eta_n \overset{\mathrm{tv}}{\to} \eta_\infty \) as \( n \to \infty \).

We say that $\mathcal{F}$ is \textit{universally bounded} if $\kappa(\mathcal{F}) \triangleq \sup\{|f(x)|:x\in S, f\in \mathcal{F}\}<\infty$. Given a transition kernel \( Q \) on \( S \), we say that \( Q \) has the \textit{\( \mathcal{F} \) inclusion property} if \( Qf \in \mathcal{F} \) whenever \( f \in \mathcal{F} \), where
\[
(Qf)(x) = \int_S f(y) Q(x, dy)
\]
for \( x \in S \).

\textbf{Remark 2.2.} Note that every transition kernel \( Q \) has the \( \mathcal{F}_e \) inclusion property, since \( |f(x)| \leq 1 \) for \( x \in S \) always implies that \( |(Qf)(x)| \leq 1 \) for \( x \in S \).

\textbf{Theorem 1.} Suppose $\mathcal{F}$ is universally bounded. For \( 1 \leq n \leq \infty \), suppose that \( P_n \) possesses a unique stationary distribution \( \pi_n \) and has the \( \mathcal{F} \) inclusion property. Assume that for each \( x \in S \) and \( m \geq 1 \):
\begin{enumerate}
    \item[(i)] \( \| P_{n,x}(X_1 \in \cdot) - P_{\infty,x}(X_1 \in \cdot) \|_{\mathcal{F}} \to 0 \) as \( n \to \infty \);
    \item[(ii)] \( \| P_{\infty,x}(X_m \in \cdot) - \pi_\infty(\cdot) \|_{\mathcal{F}} \to 0 \) as \( m \to \infty \);
    \item[(iii)] \( \| \pi_n - \pi_\infty \|_{\mathcal{F}} \to 0 \) as \( n \to \infty \).
\end{enumerate}
Then, for each \( x \in S \),
\begin{equation}
    \sup_{m \geq 0} \| P_{n,x}(X_m \in \cdot) - P_{\infty,x}(X_m \in \cdot) \|_{\mathcal{F}} \to 0 \label{eq:2.1}
\end{equation}
as $n \to \infty$.

We note that \eqref{eq:2.1} implies that if \( m_n \to \infty \), \( \| P_{n,x}(X_{m_n} \in \cdot) - \pi_\infty(\cdot) \|_{\mathcal{F}} \to 0 \), so that it implies  diagonal convergence. An interesting aspect of the theorem is that it does not assume that \( \| P_{n,x}(X_m \in \cdot) - \pi_n(\cdot) \|_{\mathcal{F}} \to 0 \) as \( m \to \infty \).

\textbf{Proof of Theorem 1.} We first note that if $f\in\mathcal{F}$, then for $m\ge 2$,
\begin{align}
  &  \big|\mathbb{E}_{n,x}f(X_m) - \mathbb{E}_{\infty,x}f(X_m)\big| \nonumber\\
  & = \bigg|\int_S P_{\infty,x}(X_{m-1} \in dy) \big[(P_{\infty} f)(y) - (P_n f)(y)\big] \nonumber\\
  & \quad + \int_S P_{\infty,x}(X_{m-1} \in dy)(P_n f)(y)
- \int_S P_{n,x}(X_{m-1} \in dy)(P_n f)(y)\bigg| \nonumber \\
& \leq \int_S P_{\infty,x}(X_{m-1} \in dy) 
\big\|P_{\infty,y}(X_{1}\in \cdot) - P_{n,y}(X_{1}, \cdot)\big\|_{\mathcal{F}} \nonumber\\
&\quad 
+ \big\|P_{\infty,x}(X_{m-1}\in \cdot) - P_{n,x}(X_{m-1}\in \cdot)\big\|_{\mathcal{F}}, \label{eq:+11}
\end{align}
where we use the fact that $P_n f\in \mathcal{F}$ for the inequality. Observe that
\[
\| P_{\infty,y}(X_1\in \cdot) - P_{n,y}(X_1\in \cdot) \|_{\mathcal{F}} \le 2 \kappa(\mathcal{F}),
\]
and hence i), in combination with the Bounded Convergence Theorem, shows that the first term on the right-hand side of \eqref{eq:+11} converges to 0. Induction in $m$ therefore establishes that i) implies that
\[
\| P_{\infty,x}(X_m\in \cdot) - P_{n,x}(X_m\in \cdot) \|_{\mathcal{F}}\rightarrow 0
\]
as $n\to\infty$. 

We next observe that for \( 1 \leq n \leq \infty \) and \( m \geq 0 \),
\begin{align}
\| P_{n,x}(X_{m+1} \in \cdot) - \pi_n(\cdot) \|_{\mathcal{F}}   =\, &  \sup_{f \in \mathcal{F}} \left| \int_S f(y) \big( P_{n,x}(X_{m+1} \in dy) - \pi_n(dy) \big) \right| \nonumber \\
  =\, & \sup_{f \in \mathcal{F}} \left| \int_S \big( P_{n}f \big)(y) \big( P_{n,x}(X_m \in dy) - \pi_n(dy) \big) \right| \nonumber \\
  =\, & \sup_{ P_nf: f\in \mathcal{F}} \left| \int_S g(y) \big( P_{n,x}(X_m \in dy) - \pi_n(dy) \big) \right| \nonumber\\
  \le\, & \sup_{g \in \mathcal{F}} \left| \int_S g(y) \big( P_{n,x}(X_m \in dy) - \pi_n(dy) \big) \right| \nonumber\\
  =\,&\| P_{n,x}(X_m \in \cdot) - \pi_n(\cdot) \|_{\mathcal{F}}, \label{eq:2.2}
\end{align}
so $(\| P_{n,x}(X_m \in \cdot) - \pi_n(\cdot) \|_{\mathcal{F}}: m\ge 0)$ is a non-increasing sequence. (We used the stationarity of \( \pi_n \) in the second equality and the \( \mathcal{F} \)-inclusion property for the inequality.)

For \( \epsilon > 0 \), property (ii) implies the existence of \(t= t(\epsilon) \) for which
\[
\| P_{\infty,x}(X_{t} \in \cdot) - \pi_\infty(\cdot) \|_{\mathcal{F}} < \epsilon / 2.
\]
Since \( \| \cdot \|_{\mathcal{F}} \) satisfies the triangle inequality,
\begin{align}
    & \| P_{n,x}(X_m \in \cdot) - P_{\infty,x}(X_m \in \cdot) \|_{\mathcal{F}}  \nonumber\\
    \le \, &  \max_{0 \leq j < t} \| P_{n,x}(X_j \in \cdot) - P_{\infty,x}(X_j \in \cdot) \|_{\mathcal{F}}
+ \sup_{j \geq t} \| P_{n,x}(X_j \in \cdot) - P_{\infty,x}(X_j \in \cdot) \|_{\mathcal{F}}\nonumber \\
\leq \,& \max_{0 \leq j < t} \| P_{n,x}(X_j \in \cdot) - P_{\infty,x}(X_j \in \cdot) \|_{\mathcal{F}}
+ \sup_{j \geq t} \| P_{n,x}(X_j \in \cdot) - \pi_n(\cdot) \|_{\mathcal{F}} \nonumber\\
&
+ \| \pi_n(\cdot)- \pi_\infty(\cdot) \|_{\mathcal{F}}
+ \sup_{j \geq t} \| P_{\infty,x}(X_j \in \cdot) - \pi_\infty(\cdot) \|_{\mathcal{F}}. \nonumber \\
\leq \, & \max_{0 \leq j < t} \| P_{n,x}(X_j \in \cdot) - P_{\infty,x}(X_j \in \cdot) \|_{\mathcal{F}}
+ \| P_{n,x}(X_t \in \cdot) - \pi_n(\cdot) \|_{\mathcal{F}} \nonumber\\
& + \| \pi_n - \pi_\infty \|_{\mathcal{F}} 
+ \| P_{\infty,x}(X_t \in \cdot) - \pi_\infty(\cdot) \|_{\mathcal{F}}\nonumber \\
=\,& \max_{0 \leq j < t} \| P_{n,x}(X_j \in \cdot) - P_{\infty,x}(X_j \in \cdot) \|_{\mathcal{F}} \nonumber \\
&+ \| (P_{n,x}(X_t \in \cdot) - P_{\infty,x}(X_t \in \cdot)) - (\pi_n(\cdot) - \pi_\infty(\cdot)) 
\nonumber\\
&
+ P_{\infty,x}(X_t\in \cdot)-\pi_\infty (\cdot) \|_{\mathcal{F}} + \| \pi_n - \pi_\infty \|_{\mathcal{F}}  + \| P_{\infty,x}(X_t \in \cdot) - \pi_\infty(\cdot) \|_{\mathcal{F}}\nonumber\\
\leq & 2\max_{0 \leq j \le t} \| P_{n,x}(X_j \in \cdot) - P_{\infty,x}(X_j \in \cdot) \|_{\mathcal{F}} 
+ 2 \| \pi_n - \pi_\infty \|_{\mathcal{F}} + 2 \|P_{\infty,x}(X_t\in\cdot) - \pi_\infty(\cdot)\|_{\mathcal{F}} \nonumber\\
\le & 
2 \max_{0 \leq j \le t} \| P_{n,x}(X_j \in \cdot) - P_{\infty,x}(X_j \in \cdot) \|_{\mathcal{F}} 
 + 2 \| \pi_n - \pi_\infty \|_{\mathcal{F}} +\epsilon.\label{eq:2.3}
\end{align}
Here, we used \eqref{eq:2.2} for the third inequality above. Hence, \eqref{eq:2.3} implies that:
\begin{align*}
    &\sup_{m \geq 0} \| P_{n,x}(X_m \in \cdot) - P_{\infty,x}(X_m \in \cdot) \|_{\mathcal{F}} \\
&\leq 2\max_{0 \leq j \le t} \| P_{n,x}(X_j \in \cdot) - P_{\infty,x}(X_j \in \cdot) \|_{\mathcal{F}} 
+ 2 \| \pi_n - \pi_\infty \|_{\mathcal{F}} + \epsilon.
\end{align*}
We now send \( n \to \infty \) using (i) and (iii), and then send \( \epsilon \downarrow 0 \), yielding the desired conclusion. $\square$

We now apply this theorem to obtain diagonal convergence when we have convergence in total variation.

\textbf{Proposition 1.} 

Assume  that for \( 1 \leq n \leq \infty \), \( P_n \) has a unique stationary distribution \( \pi_n \), and suppose that for each \( x \in S \),
\begin{enumerate}
    \item[(i)] \( P_n(x, \cdot) \overset{\mathrm{tv}}{\to} P_\infty(x, \cdot) \) as \( n \to \infty \);
    \item[(ii)] \( P_{\infty,x}(X_m \in \cdot) \overset{\mathrm{tv}}{\to} \pi_\infty(\cdot) \) as \( m \to \infty \);
        \item[(iii)] \( \pi_n \overset{\mathrm{tv}}{\to} \pi_\infty \) as \( n \to \infty \).
\end{enumerate}
Then, for \( x \in S \),
\[
\sup_{m \geq 0} \| P_{n,x}(X_m \in \cdot) - P_{\infty,x}(X_m \in \cdot) \|_e \to 0 
\]
as $n \to \infty.$

\textbf{Proof.} We apply Theorem 1 with \( \mathcal{F} = \mathcal{F}_e \) and note that $\mathcal{F}_e$ is universally bounded. $\Box$

\textbf{Corollary 1.} Suppose that $S$ is finite or countably infinite. For $1 \leq n \leq \infty$, assume that the transition matrix $P_n = (P_n(x,y):{x,y \in S})$ has a unique stationary distribution $\pi_n = (\pi_n(x):{x \in S})$. If, for each $x, y \in S$,
\begin{itemize}
    \item[(i)] $P_n(x,y) \to P_\infty(x,y)$ as $n \to \infty$;
    \item[(ii)] $P_\infty^m(x,y) \to \pi_\infty(y)$ as $m \to \infty$;
    \item[(iii)] $\pi_n(y) \to \pi_\infty(y)$ as $n \to \infty$,
\end{itemize}
then
\[
\sup_{m \geq 0} \sum_{y \in S} \left| P_n^m(x,y) - P_\infty^m(x,y) \right| \to 0 \quad \text{as } n \to \infty.
\]

The corollary follows immediately from Scheffé's lemma (see, for example, p.246 of \cite{billingsley1968convergence}), since
it shows that point-wise convergence of $\nu$-densities implies $L^1$ convergence. So, when $S$ is discrete, the weak interchange property implies diagonal convergence.

\textbf{Remark 2.3.} In the truncation context for countably infinite chains, it is standard that the $n$'th truncation corresponds to a finite state space $S_n$, where $\phi \neq S_1 \subseteq S_2 \subseteq ... ,$ and $\cup_{n=1}^{\infty} S_n = S$, and $(P_n(x,y): x,y \in S_n)$ is a transition matrix with a single closed communicating class (so that $(\pi_n(x): x \in S_n)$ exists uniquely). To apply the corollary to the truncation setting, we extend $(P_n(x,y): x,y \in S_n)$ and $(\pi_n(x): x \in S_n)$ to $S$ by choosing an arbitrary $z \in S_1$ and putting $P_n(x,z) = 1$ for $x \notin S_n$ and $\pi_n(x) = 0$ for $x \notin S_n$. This extension makes $S_n$ absorbing for $P_n$. We conclude that whenever a truncation algorithm exhibits stationary distribution convergence, then $P_\infty^m(x,y)$ may be approximated by $P_n^m(x,y)$ across all time values at $m$ for $n$ sufficiently large and $x,y$ fixed, 
so that the truncation approximates well all the marginals of the Markov chain under $P_\infty$ at all time scales. (We note that without such a result, one would need to worry about the range of $m$ over which good approximations hold.)

We now extend Corollary 1 to the more general weighted total variation norm $\|\cdot\|_w$ when $S$ is discrete.

\textbf{Corollary 2.} Suppose that $S$ is finite or countably infinite. If the conditions of Corollary 1 hold, $P_\infty$ is irreducible, and 
\begin{equation}
    \sum_{x \in S} w(x) \pi_n(x) \to \sum_{x \in S} w(x) \pi_\infty(x) \label{eq:2.6}
\end{equation}
as $n \to \infty$, then
\[
\sup_{m \geq 0} \sum_{y \in S} w(y) \left| P_n^m(x,y) - P_\infty^m(x,y) \right| \to 0
\]
as $n \to \infty.$

\textbf{Proof.} We observe that the stationarity of $\pi_n$ implies that
\[
\sum_{x \in S} \pi_n(x) P_n^m(x,y) = \pi_n(y).
\]
for $y \in S$ and $m, n \geq 1$. So,
\[
\pi_n(x) P_n^m(x,y) \leq \pi_n(y).
\]
The irreducibility of $P_\infty$ implies that $\pi_n(x) > 0$ for $n\ge n_0$ with $n_0$ sufficiently large. Hence,
\begin{equation}
    P_n^m(x,y) \leq \frac{\pi_n(y)}{\pi_n(x)} \label{+12}
\end{equation}
for $x, y \in S$ and $n\ge n_0$. Since $\pi_n \overset{tv}{\to} \pi_\infty$ and \eqref{eq:2.6} holds, $w$ is uniformly integrable with respect to $(\pi_n:n \geq 1)$ (see \cite{billingsley1968convergence}), so that 
\[
\sup_{n_0 \leq j \leq \infty} \sum_{y} w(y) \pi_j(y) \mathbf{I}(w(y) \geq b) \to 0
\]
as $r \to \infty$. Hence, in view of \eqref{+12}, 
\begin{align*}
    \sum_{y \in S} w(y) \left| P_n^m(x,y) - P_\infty^m(x,y) \right| \le &  b\max_{j \geq 0} \| P_{n, x}(X_j(\cdot)) - P_{\infty, x}(X_j(\cdot)) \|_e \\
   & + 2 \sup_{n_0\leq j\leq \infty} \sum_{y \in S} w(y) \pi_j(y) \mathbf{I}(w(y) > b) / \min_{n_0\leq j\leq \infty} \pi_j(x).
\end{align*}
Sending first $n \to \infty$  (and applying Corollory 1) and then $b \to \infty$ yields the corollary. $\Box$

We conclude this section with a version of Proposition 1 that involves weak convergence rather than total variation convergence. Assume now that $S$ is a metric space with metric $d$. For $r >0$, let $\mathrm{Lip}_b(r)$ be the family of functions $f$ such that $\|f\| \leq e$ and 
\begin{equation}
    |f(x) - f(y)| \leq r \,d(x, y) \label{eq:+13}
\end{equation}
for $x, y \in S$, so that $f$ has Lipschitz constant at most $r$. Note that $\mathrm{Lip}_b(r)$ is a universally bounded family of functions. We now assume that for $1 \leq n \leq \infty$, $P_n$ has the $\mathrm{Lip}_b(1)$ inclusion property. In particular, this holds if the Markov chain $X$ can be represented as a random iterated function system under $P_n$, in which, conditional on $X_0 = x$,
\[
X_m = (\varphi_n(m) \circ \varphi_{n}(m-1) \circ \dots \circ \varphi_n(1))(x),
\]
where $(\varphi_n(i) : i \geq 1)$ is a sequence of independent and identically distributed (i.i.d.) random functions mapping $S$ into $S$ such that 
\[
\E d(\varphi_n(1)(x), \varphi_n(1)(y)) \leq r\,d(x, y)
\]
for $x, y \in S$, where $0\le r\le 1$. To verify the inclusion property, note that for $f \in \mathrm{Lip}_b(1)$,
\begin{align*}
    &|(P_n f)(x) - (P_n f)(y)| \\
    \leq & \mathbb{E}|f(\varphi_n(1)(x)) - f(\varphi_n(1)(y))|\\
    \leq & \mathbb{E} |d(\varphi_n(1)(x), \varphi_n(1)(y))|\\
    \leq & r\,d(x, y).
\end{align*}
for $x,y\in S$, so that $P_n f \in \mathrm{Lip}_b(r) \subseteq \mathrm{Lip}_b(1)$. 

\textbf{Remark 2.4.} Suppose that $\mathrm{Lip}(r)$ is the family of functions satisfying \eqref{eq:+13} (but with no requirement that $|f|\le e$). It is worth noting that $\|\eta_1-\eta_2\|_{\mathrm{Lip}(1)}$ is then exactly the Wasserstein 1-distance $W_1(\eta_1,\eta_2)$. 
    
Let $\Rightarrow$ denote weak convergence in $S$.

\textbf{Proposition 2.} Suppose that $S$ is a separable metric space with metric $d$, and assume that for $1 \leq n \leq \infty$, $P_n$ has the $\mathrm{Lip}_b(1)$ inclusion property and possesses a unique stationary 
distribution $\pi_n$. Assume further that for each $x \in S$,

\begin{itemize}
    \item[(i)] $P_n(x, \cdot) \Rightarrow P_\infty(x, \cdot)$ as $n \to \infty$;
    \item[(ii)] $P_{\infty, x}(X_m\in \cdot) \Rightarrow \pi_\infty(\cdot)$ as $m \to \infty$;
    \item[(iii)] $\pi_n \Rightarrow \pi_\infty$ as $n \to \infty$.
\end{itemize}
Then,
\begin{equation}
    \sup_{m \geq 0} \| P_{n, x}(X_m\in\cdot) - P_{\infty, x}(X_m\in \cdot) \|_{\text{Lip}_b(1)} \to 0 \quad \text{as } n \to \infty.  \label{eq:2.7}
\end{equation}

\textit{Proof.} We will apply Theorem 1 with $\mathcal{F} = \text{Lip}_b(1)$. All that is then needed to finish the proof is to establish that weak convergence implies convergence in $\|\cdot\|_{\text{Lip}_b(1)}$. To this end, if $\eta_n\Rightarrow \eta_\infty$ as $n\rightarrow\infty$, the Skorohod Representation Theorem ensures that there exists $(\Lambda_n:1\le n\le \infty)$ such that $\Lambda_n$ has distribution $\eta_n$ for $1\le n\le \infty$ and $\Lambda_n\rightarrow\Lambda_\infty$ a.s. as $n\rightarrow\infty$. Then, for $f\in \text{Lip}_b(1)$, 
\[
|\int_S (\eta_n(dx) - \eta_\infty(dx))f(x)| = |\E(f(\Lambda_n) - f(\Lambda_\infty))| \le |\E(d(\Lambda_n,\Lambda_\infty)\land 2)|
\]
uniformly in $f$ (where $a\land b \triangleq \min(a,b)$. We used the boundedness and Lipschitzness of $f$ for the inequality. Use of the Bounded Convergence Theorem then proves that $\|\eta_n-\eta_\infty\|_{\text{Lip}_b(1)} \to 0$ as $n\to\infty$.  $\Box$

\textbf{Example 2.1.} Suppose that $(X_j:j\geq 0)$ is a Markov chain describing the delay sequence for the single-server $GI/G/1$ queue, so that under $P_n$,
\[
(P_nf)(x) = \mathbb{E}f([x + Z_n]^+)
\]
for some random variable (r.v.) $Z_n$, where $[y]^{+} \triangleq \max(y, 0)$.
If $f \in \text{Lip}_b(1)$, then
\[
|(P_nf)(x) - (P_nf)(y)| \leq E[|[x + Z_n]^+ - [y + Z_n]^+|]\leq |x - y|,
\]
so $P_nf \in \text{Lip}_b(1)$, and hence $P_n$ has the $\text{Lip}_b(1)$ contraction property. If $Z_n \nearrow Z_\infty$ as $n \to \infty$ with $\mathbb{E} Z_\infty < 0$, then it is well known that $\pi_n$ exists for $1 \leq n \leq \infty$ (see \cite{lindley1952theory} and \cite{kiefer1955theory}), and conditions (i) and (ii) of Proposition 2 are easily verified.

For condition (iii), let $((Z_n(i): n \geq 1): i \geq 1)$ be an i.i.d. sequence of copies of $\{Z_n : n \geq 1\}$ and put $S_n(k) = \sum_{j=1}^k Z_n(j)$. Then,
\[
\pi_n(\cdot) = P(\max_{k \geq 0} S_n(k) \in\cdot),
\]
(see \cite{kiefer1955theory}), and
\[
\max_{k \geq 0} S_n(k) = \max_{0 \leq k \leq L_n} S_n(k),
\]
where $L_n = \max \{k \geq 0 : S_n(k) > 0 \} < \infty$. Since $Z_n(j) \leq Z_\infty (j)$, $L_n \leq L_\infty$. So,
\[
\max_{k \geq 0} S_n(k) = \max_{0 \leq k \leq L_\infty} S_n(k) \nearrow \max_{0 \leq k \leq L_\infty} S_\infty(k),
\]
so that $\pi_n \Rightarrow \pi_\infty$ as $n \to \infty$. Consequently, \eqref{eq:2.7} holds for this model, so that diagonal convergence follows in particular. This result both complements earlier findings for the $M/G/1$ queue (see \cite{breuer2008continuity}) and strengthens the conclusions reached there (via application of Proposition 2).

\textbf{Example 2.2.} Suppose that $(X_j : j \geq 0)$ is a \textit{contractive} Markov chain on $S \subseteq \mathbb{R}^d$, for which $P_n$ can be represented as
\begin{equation}
    P_n(x, \cdot) = P(\varphi_n(x)\in \cdot), \label{eq:+15}
\end{equation}
where $(\varphi_n : n \geq 1)$ is a sequence of continuous random mappings $\varphi_n : S \to S$ defined on a common probability space for which:
\begin{itemize}
    \item[(a)] $\mathbb{E} | \varphi_n(x) - \varphi_n(y)\| | \leq r \|x - y\|$;
    \item[(b)] $\mathbb{E} |\varphi_n(x) - x| \leq \infty$;
    \item[(c)] $\varphi_n\to \varphi_\infty$  uniformly on compact sets as $n \to \infty$;
    \item[(d)] $\mathbb{E} |\varphi_n(x) - \varphi_\infty(x)|  \to 0$ for $n \ge 1$
\end{itemize}
for all $x, y \in S$, where $r < 1$.

As noted earlier, (b) implies that $P_n$ has the $\text{Lip}_b(1)$ inclusion property for $1\le n \le \infty$. Hence, Proposition 2 is applicable, once we verify (i)–(iii).

Assumption (c) ensures that (i) holds.  For (ii), let $(\varphi_n(k):n\ge 1) : k \geq 1\}$ be a sequence of independent copies of $(\varphi_n : n \geq 1)$.

For $x \in S$, let
\[
\beta_n(k, x) = \big(\varphi_n(1) \circ \varphi_n(2) \circ \cdots \circ \varphi_n(k)\big)(x),
\]
and note that
\[
P_{n,x}(X_k \in \cdot) = P(\beta_n(k, x) \in \cdot).
\]

Under our contraction condition (a), it is known that there exists $\beta_n(\infty)$ such that
\begin{equation}
\mathbb{E}|\beta_n(\infty) - \beta_n(k, x)| \leq \frac{r^k}{1 - r} \mathbb{E}|\varphi_n(x) - x|, \label{eq:+18}
\end{equation}
so that
\[
P_{n,x}(X_k \in \cdot) \Rightarrow \pi_n(\cdot) = P(\beta_n(\infty) \in \cdot)
\]
as $k\rightarrow\infty$, validating ii) for $n=\infty$; see, for example, Section 5 of \cite{diaconis1999iterated} for the bound \eqref{eq:+18}.

For iii), we start by noting that c) ensures that $$
\beta_n(k, x) = \varphi(1,\varphi_n(k,x))\to \varphi_\infty(1,\varphi_\infty(k,x)) = \beta_\infty(k,x)$$
as $n \to \infty$. Similarly, we conclude that $\beta_n(k, x) \to \beta_{\infty}(k,x)$ as $n \to \infty$ for $k \ge 3$. For $k \ge 1$ and $\epsilon > 0$, we now write
\begin{align}
    & P\big(|\beta_n(\infty) - \beta_n(k, x)| > \epsilon\big)  \nonumber\\
    &\leq 
P\big(|\beta_n(\infty) - \beta_n(k, x)| + |\beta_{\infty}(\infty) - \beta_\infty(k, x)| > \epsilon/2\big)
+ P\big(|\beta_n(k, x) - \beta_{\infty}(k, x)| > \epsilon / 2\big). \label{eq:+3} 
\end{align}
For the first term on the right-hand side, we use \eqref{eq:+18} and Markov's inequality to obtain:
\begin{align}
  &  P\big(|\beta_n(\infty) - \beta_n(k, x)| + |\beta_{\infty}(k,x) - \beta_\infty(\infty)| > \epsilon / 2\big) \nonumber \\
  \le & \frac{2}{\epsilon} \mathbb{E}|\beta_n(\infty) - \beta_n(k, x)| + |\beta_{\infty}(\infty) - \beta_\infty(k, x)| \nonumber\\
  \le &\frac{2}{\epsilon} \frac{r^k}{1-r} (\mathbb{E}|\varphi_n(x) - x| + |\varphi_\infty(x) - x|). \label{eq:+44} 
\end{align}
We now let $n \to \infty$ in \eqref{eq:+3}, using the upper bound \eqref{eq:+44} and d), thereby yielding 
\[
\overline{\lim}_{n \to \infty} P\big(|\beta_n(\infty) - \beta_\infty(\infty)| > \epsilon\big) 
\leq \frac{4 r^k}{\epsilon(1 - r)} \mathbb{E}|\varphi_\infty(x) - x|.
\]
Since $k$ was arbitrary, we may now send $k \to \infty$, use (b), and conclude that $\beta_n(\infty) \overset{p}{\to} \beta_{\infty}(\infty) \text{ as } n \to \infty,$ yielding (iii).

\section{Convergence of First Transition Expectations}
In this section, we focus exclusively on the case where $S$ is finite or countably infinite. Our goal is to show that when a sequence $(P_n: 1 \leq n \leq \infty)$ is such that $P_n \to P_\infty$ as $n \to \infty$, and $\pi_n \to \pi_\infty$, then it is guaranteed that expected hitting times, expected infinite horizon discounted rewards, and many other associated expectations computed for $P_n$ converge to those associated with $P_\infty$. As noted in the Introduction, stationary distribution convergence then allows us to approximate these $P_n$ expectations by those associated with the limiting model $P_\infty$ or, equivalently, in computing $P_\infty$ expectations, one can use the $P_n$ expectations as approximations (as in the truncation setting; see Remark 2.3).

Given $\emptyset \neq C \subseteq S$, let 
\[
T = \inf\{n \geq 0 : X_n \in C^c\}
\]
be the first hitting time of $C^c$. (Note that $C^c$ may be empty, in which case $T = \infty$.) For $r : S \to \mathbb{R}_+$ and $\alpha : S \rightarrow [0, \infty)$, let
\begin{equation}
    u_n^*(x) = \mathbb{E}_{n, x} \left[ \sum_{j=0}^{T} \exp\left( - \sum_{k=0}^{j-1} \alpha(X_k) \right) r(X_j) \right] \label{eq:3.1}
\end{equation}
and note that $u_n^* = \left( u_n^*(x) : x \in S \right)$ is the minimal non-negative solution of
\begin{equation}
    u_n(x) = r(x) + \sum_{y \in C^c} G_n(x, y) r(y) + \sum_{y \in C} G_n(x, y) u_n(y)  \label{eq:3.2}
\end{equation}
for \(x \in C\), where \(G_n = \left(G_n(x, y) : x, y \in S \right)\) 
has entries given by $$G_n(x, y)= \exp\left(-\alpha(x)\right)P_n(x, y)$$ for $x, y \in S$. We say that \(u_n^*(x)\) is a first transition expectation, since the linear system \eqref{eq:3.2} can be derived from \eqref{eq:3.1} by conditioning on the first transition state value \(X_1\).

If we set \(\alpha(y) = 0\) for \(y \in S\) and \(r = e\), then \(u_n^*(x)\) is the mean hitting time of \(C_n^c\) starting from \(x\). On the other hand, if \(\alpha(y) = \alpha > 0\) and \(C^c = \emptyset\), then \(u_n^*(x)\) is the expected infinite horizon discounted reward starting from \(x\).
If \( r(y) = 0 \) for \( y \in C \), then \( u_n^*(x) \) is the \(\alpha\)-discounted reward collected when the Markov chain hits \( C^c \), where \(\alpha(\cdot)\) represents a state-dependent risk-free rate (as occurs in many finance settings, in which \( C^c \) is the exercise region for an option).

\textbf{Theorem 2.} Suppose that \( S \) is finite or countably infinite, and that $r \in \mathcal{F}_e$. Assume that $P_\infty = \left( P_\infty(x, y) : x, y \in S \right)$ is irreducible and that \( P_n = \left( P_n(x, y) : x, y \in S \right) \) has a unique stationary distribution \(\pi_n\) for $1 \leq n \leq \infty$. If, for each \( x, y \in S \),
\begin{itemize}
    \item[(i)] \( P_n(x, y) \to P_\infty(x, y) \) as \( n \to \infty \);
    \item[(ii)] \( \pi_n(x) \to \pi_\infty(x) \) as \( n \to \infty \),
\end{itemize}
then, for each \( x \in S \), \( u_n^*(x) \to u_\infty^*(x) \) as \( n \to \infty \).

\textbf{Proof.} We first note that condition (i) ensures that $P_n(x, \cdot) \overset{\text{tv}}{\to} P_\infty(x, \cdot)$ as $n \to \infty.$ The proof of Theorem 1 can be easily modified to ensure that for \( m \geq 1 \),
\[
\| P_{n, x} \left( (X_0, X_1, \ldots, X_m) \in \cdot \right) - P_{\infty, x} \left( (X_0, X_1, \ldots, X_m) \in \cdot \right) \|_e \to 0
\]
as \( n \to \infty \), from which it follows easily that 
\[
P_{n, x} \left( \sum_{j=0}^{T} \exp \left( -\sum_{k=0}^{j-1} \alpha(X_k) \right) r(X_j) \in \cdot \right) 
\]

\begin{equation}
   \implies  P_{\infty, x} \left( \sum_{j=0}^{T} \exp \left( -\sum_{k=0}^{j-1} \alpha(X_k) \right) r(X_j) \in \cdot \right) \label{eq:3.3}
\end{equation}
as \( n \to \infty \).

Next, if \( x \in C^c \), then \( u_n^*(x) = r(x) \) for $1 \leq n \leq \infty$, so the result is trivial. If \( x \in C \), let \( \tau(x) = \inf \{ n \geq 1 : X_n = x \} \) be the first entry time into \( x \), and note that the strong Markov property implies that 

\[
u_n^*(x) = \mathbb{E}_{n, x} \left[ \sum_{j=0}^{T \land (\tau(x)-1)} \exp\left(-\sum_{k=0}^{j-1} \alpha(X_k)\right) r(X_j) \right]
+ 
\]
\[
\mathbb{E}_{n, x} \left[ \exp\left(-\sum_{k=0}^{\tau(x)-1} \alpha(X_k)\right) \mathbf{1}\{T > \tau(x)\} u_n^*(x) \right].
\]
Note that since \( r \in \mathcal{F}_e \) and \( \alpha(\cdot) \) is non-negative,
\begin{equation}
    \sum_{j=0}^{T \land (\tau(x)-1)} \exp\left(-\sum_{k=0}^{j-1} \alpha(X_k)\right) r(X_j) \leq \tau(x), \label{eq:3.4}
\end{equation}
and
\begin{equation}
    \exp\left(-\sum_{k=0}^{\tau(x)-1} \alpha(X_k)\right) \mathbf{1}\{T > \tau(x)\} \leq 1. \label{eq:3.5}
\end{equation}
Because $\pi_n(x) = \left( \mathbb{E}_{n, x} \tau(x) \right)^{-1} \to \pi_\infty(x) = \left( \mathbb{E}_{\infty, x} \tau(x) \right)^{-1} > 0 $ (as a consequence of (ii) and irreducibility), \eqref{eq:3.3} implies that the left-hand sides of \eqref{eq:3.4} and \eqref{eq:3.5} are uniformly integrable over the family of probabilities $\left(P_{n,x} : n \geq n_0 \right)$ (for \( n_0 \) sufficiently large). Hence,
\begin{equation}
    \mathbb{E}_{n,x} \left[ \sum_{k=0}^{T \land (\tau(x)-1)} \exp\left(-\sum_{j=0}^{k-1} \alpha(X_j)\right) r(X_k) \right] 
\to \mathbb{E}_{\infty,x} \left[ \sum_{k=0}^{T\land (\tau(x)-1)} \exp\left(-\sum_{j=0}^{k-1} \alpha(X_k)\right) r(X_j) \right].  \label{eq:3.6}
\end{equation}
and
\begin{equation}
    \mathbb{E}_{n,x} \left[ \exp\left(-\sum_{j=0}^{\tau(x)-1} \alpha(X_k)\right) \mathbf{1}\{T > \tau(x)\} \right] 
\to \mathbb{E}_{\infty,x} \left[ \exp\left(-\sum_{j=0}^{\tau(x)-1} \alpha(X_j)\right) \mathbf{1}\{T > \tau(x)\} \right]. \label{eq:3.7}
\end{equation}
as $n\rightarrow\infty$. If \( u_n^*(x) < \infty \), then
\begin{equation}
    u_n^*(x) = \frac{\mathbb{E}_{n,x} \left[ \sum_{k=0}^{T \land (\tau(x)-1)} \exp\left(-\sum_{j=0}^{k-1} \alpha(X_k)\right) r(X_j) \right] }
{1 - \mathbb{E}_{n,x} \left[ \exp\left(-\sum_{j=0}^{\tau(x)-1} \alpha(X_j)\right) \mathbf{1}\{T > \tau(x)\} \right]}. \label{eq:3.8}
\end{equation}
It follows from \eqref{eq:3.6} and \eqref{eq:3.7} that if \( u_\infty^*(x) < \infty \), then \( u_n^*(x) \to u_\infty^*(x) \). On the other hand, if \( u_\infty^*(x) = \infty \), then 
\[
\mathbb{E}_{\infty,x} \left[ \exp\left(-\sum_{j=0}^{\tau(x)-1} \alpha(X_j)\right) \mathbf{1}\{T > \tau(x)\} \right] = 1,
\]
and
\[
\mathbb{E}_{\infty,x} \left[ \sum_{j=0}^{T \land (\tau(x)-1)} \exp\left(-\sum_{k=0}^{j-1} \alpha(X_k)\right) r(X_j) \right] > 0,
\]
so that
\eqref{eq:3.6}, \eqref{eq:3.7}, and \eqref{eq:3.8} again imply that 
\[
u_n^*(x) \to \infty = u_\infty^*(x). 
\]
$\Box$

Hence, stationary distribution convergence automatically implies convergence of all expectations covered by Theorem 2. An implication is that in developing convergent truncation schemes, it is enough to establish stationary distribution convergence in order to establish ``universal convergence'' across all the expectations considered above.

We conclude this section with a modest extension to Theorem 2, covering the case where \( r \) is unbounded.

\textbf{Proposition 3.} Suppose \( P \) is irreducible, that \( P_n \) has a stationary distribution \(\pi_n\) for \( 1\le n \le \infty  \), and that (i) and (ii) of Theorem 2 hold. If, in addition,
\begin{equation}
    \sum_{x \in S} \pi_n(x) w(x) \to \sum_{x \in S} \pi_\infty(x) w(x)\label{eq:3.9}
\end{equation}
as $n\to\infty$, then \( u_n^*(x) \to u_\infty^*(x) \) as \( n \to \infty \) for each \( r \in \mathcal{F}_w \).

\textbf{Proof.} Since \( \pi_n(x) \to \pi_\infty(x) > 0 \) as \( n \to \infty \) as a consequence of (ii) and irreducibility, \( x \) must lie in \( P_n \)'s positive recurrent closed communicating class for $n$ sufficiently large. Applying regenerative process theory, we find that 
\[
\sum_{y \in S} \pi_n(y) w(y) = \frac{\mathbb{E}_{n,x} \left[ \sum_{j=0}^{\tau(x)-1} w(X_j) \right]}{\mathbb{E}_{n,x} \tau(x)};
\]
for \( 1 \leq n \leq \infty \), see \cite{asmussen2003applied}. Hence, \eqref{eq:3.9} implies that 
\[
\frac{\mathbb{E}_{n,x} \left[ \sum_{j=0}^{\tau(x)-1} w(X_j) \right]}{\mathbb{E}_{n,x} \tau(x)} = \sum_{y \in S} \pi_n(y) w(y)
\]
\[
\rightarrow \sum_{y \in S} \pi_\infty(y) w(y) = \frac{\mathbb{E}_{\infty,x} \left[ \sum_{j=0}^{\tau(x)-1} w(X_j) \right]}{\mathbb{E}_{\infty,x} \tau(x)}.
\]
as $n\to\infty$. Since \( \mathbb{E}_{n,x} \tau(x) \to \mathbb{E}_{\infty,x} \tau(x) \) as \( n \to \infty \), we find that 
\[
\mathbb{E}_{n,x} \left[ \sum_{j=0}^{\tau(x)-1} w(X_j) \right] \to \mathbb{E}_{\infty,x} \left[ \sum_{j=0}^{\tau(x)-1} w(X_j) \right]
\]
as $n\to\infty$. In view of \eqref{eq:3.3} with \( \alpha(\cdot) = 0 \), \( T = \tau(x) \), and \( r = w \), we may therefore conclude that $\sum_{j=0}^{\tau(x)-1} w(X_j) $ is uniformly integrable under \( ( P_{n,x} : 1 \leq n \leq \infty ) \). If \( r \in \mathcal{F}_w \), then 
\[
\sum_{k=0}^{T\land (\tau(x)-1)} \exp\left(-\sum_{j=0}^{k-1} \alpha(X_j)\right) r(X_k) \leq \sum_{k=0}^{\tau(x)-1} w(X_k),
\]
so $\sum_{k=0}^{T\land(\tau(x)-1)} \exp\left(-\sum_{j=0}^{k-1} \alpha(X_j)\right) r(X_k)$ is also uniformly integrable under \( ( P_{n,x} :1 \leq n \leq \infty ) \). Utilizing \eqref{eq:3.3} with \( T \) replaced by \( T \wedge (\tau(x) - 1) \), we conclude that 
\[
\mathbb{E}_{n,x} \left[ \sum_{k=0}^{T \wedge (\tau(x)-1)} \exp\left(-\sum_{j=0}^{k-1} \alpha(X_j)\right) r(X_k) \right] 
\to 
\mathbb{E}_{\infty,x} \left[ \sum_{k=0}^{T \wedge (\tau(x)-1)} \exp\left(-\sum_{j=0}^{k-1} \alpha(X_j)\right) r(X_k) \right].
\]
as $n\to\infty$. We now apply \eqref{eq:3.8}, after which the rest of the argument follows identically as in the proof of Theorem 2.
$\Box$

\section{A Counter-example}
As shown in Corollary 1, the weak interchange property automatically implies diagonal convergence when \( S \) is countably infinite. In this section, we show that this can fail to be true when \( S \) is a continuous state space.

Let \( S = [0, 1] \). For \( 1 \leq n < \infty \), we put 
\[
P_n(x, \cdot) = \delta_{x/2}(\cdot)
\]
for \( x > 2^{-n} \). For \( 2^{-n-1} < x \leq 2^{-n} \), set
\[
P_n(x, \cdot) = \delta_{1}(\cdot),
\]
and for \( 0 \leq x \leq 2^{-n-1} \),
\[
P_n(x, \cdot) = \frac{1}{n+1} \sum_{j=0}^{n} \delta_{2^{-j}}(\cdot).
\]
Also, for \( n = \infty \), put \( P_\infty(x, \cdot) = \delta_{x/2}(\cdot) \) for \(0 \leq  x \leq 1 \), and note that \( \pi_\infty(\cdot) = \delta_{0}(\cdot) \) is the unique stationary distribution of \( P_\infty \).

Observe that for each \( x > 0 \),
once \( n \) is large enough that \( 2^{-n} < x \), the Markov chain starting from \( x \) under \( P_n \) takes 
\[
m_n(x) \triangleq n - \lfloor \log_2 x \rfloor
\]
steps to move into the interval \((2^{-n-1}, 2^{-n}]\), after which the chain jumps to state \( 1 \). So,
\begin{equation}
    P_{n,x} \left( X_{m_n(x)+1} \in \cdot \right) = \delta_1(\cdot). \label{eq:4.1}
\end{equation}
Once the chain hits state \( 1 \), the chain visits the states \( 2^{-1}, 2^{-2}, \ldots, 2^{-n} \) and repeats this visitation pattern from that point onwards.

On the other hand, if \( x \in (2^{-n-1}, 2^{-n}] \), the chain immediately jumps to \( 1 \), after which the chain visits \( 2^{-1}, 2^{-2}, \ldots, 2^{-n} \) and then indefinitely repeats this pattern. If \( x \in [0, 2^{-n-1}] \), the chain chooses one of the states \( 1, 2^{-1}, \ldots, 2^{-n} \) uniformly at random, after which the chain enters the periodic ``orbit" \( 1, 2^{-1}, \ldots, 2^{-n} \) and repeats that pattern forever. It follows that
\[
\frac{1}{m} \sum_{j=0}^{m-1} P_{n, x}(X_j \in \cdot) \overset{\text{tv}}{\to} \pi_n(\cdot) \triangleq \frac{1}{n+1} \sum_{j=0}^{n} \delta_{2^{-j}}(\cdot).
\]
as \( n \to \infty \). (In fact, \( X \) is a periodic positive recurrent Harris chain under \( P_n \)).

Note that
\begin{equation}
    P_n(x, \cdot) \overset{\text{tv}}{\rightarrow} P_\infty(x, \cdot) \label{eq:4.2}
\end{equation}
for \( 0 < x \leq 1 \) as \( n \to \infty \), whereas
\begin{equation}
    P_n(0, \cdot) \Rightarrow P_\infty(0, \cdot) \label{eq:4.3}
\end{equation}
as \( n \to \infty \). Also, for \( x \in S \),
\begin{equation}
    P_{\infty,x}(X_m \in \cdot) \Rightarrow \pi_\infty(\cdot) \label{eq:4.4}
\end{equation}
as \( m \to \infty \), and
\begin{equation}
    \pi_n \Rightarrow \pi_\infty. \label{eq:4.5}
\end{equation}
 as \( n \to \infty \). Consequently, \eqref{eq:4.2}-\eqref{eq:4.5} imply that the Markov chain constructed here satisfies conditions (i)-(iii) of Proposition 2 (but, critically, this example does not have the $Lip_b(1)$ inclusion property), so it has the weak interchange property. However, \eqref{eq:4.1} shows that for \( x > 0 \),
\[
P_{n,x} \left( X_{m_n(x)+1} \in \cdot \right) = \delta_1(\cdot) \rlap{\(\quad\not\)}\implies \delta_0(\cdot) = \pi_\infty(\cdot),
\]
and hence diagonal convergence fails in this example.

\textbf{Remark.} If we wish to modify this example so that \( P_n \) is Feller for \( 1 \leq n < \infty \) (so that \( P_n(x_k, \cdot) \Rightarrow P_n(x_\infty, \cdot) \) whenever \( x_k \to x_\infty \)), this is easy to do. In particular, we can modify the transition probabilities in \( (2^{-n+1}, 2^{-n}) \) and \( (0, 2^{-n-1}) \) to be mixtures of the endpoint transition probabilities, where the mixture weights are chosen to make the transition probabilities weakly continuous.

We conclude this section by noting that the theory of Section 2 depends critically on the fact that
\[
P_{\infty,x}(X_m \in \cdot) \Rightarrow \pi_\infty(\cdot)
\]
as \( m \to \infty \). This typically holds only when the chain is suitably aperiodic under \( P_\infty \). If \( P_\infty \) induces a periodic Markov chain, one can expect the result
\[
\sup_{m \geq 0} \| P_{n,x}(X_m \in \cdot) - P_{\infty, x}(X_m \in \cdot) \| \to 0
\]
to hold only when \( P_n \) matches the periodic structure of \( P_\infty \). But
\begin{equation}
    P_n(x, \cdot) \Rightarrow P_\infty(x, \cdot) \label{eq:4.6}
\end{equation}
is too weak a condition to assure that the periodic structure of \( P_n \) matches that of \( P_\infty \). (For example, when \( S \) is countably infinite, we can achieve \eqref{eq:4.6} and yet include a state \( z \in S \) for which \( P_n(z, z) = 1/n \), so that \( P_n \) is aperiodic for \( 1 \leq n < \infty \), regardless of the periodicity of \( P_\infty \).)

\section{Extension to Markov Jump Processes}

In this section, we briefly describe how our theory extends to Markov jump processes $X = (X(t) : t \geq 0)$ with a countably infinite state space $\mathcal{S}$. Here, $P_n = (P_n(x,y) : x,y \in \mathcal{S})$ for $1\le n\le\infty$ is replaced by the rate matrix 
$Q_n = (Q_n(x,y) : x,y \in \mathcal{S})$, where $Q_n$ satisfies the requirements:
\begin{enumerate}
    \item[a)] $Q_n(x,y) \geq 0 \quad \text{for } x \neq y, x,y \in \mathcal{S}$
    \item[b)] $\sum_{y\neq x} Q_n(x,y) <\infty$
    \item[c)] $\sum_{y} Q_n(x,y) = 0$
\end{enumerate}
for $x,y\in\mathcal{S}$. We further require throughout this section that $Q_n$ describes a non-explosive jump process; see Chapter II of \cite{asmussen2003applied} for the definition and sufficient conditions. As in Section 2, we let $P_{n,x}(\cdot)$ be the probability on the path-space of $X$ under which $X$ starts in state $x \in \mathcal{S}$ and evolves according to $Q_n$.  

Here is our extension of Corollary 1 to Markov jump processes.

\textbf{Proposition 4.} Suppose that for $1\le n \le \infty$, $Q_n$ has a unique stationary distribution $\pi_n$.  If:
\begin{enumerate}
    \item[i)] $Q_n(x, y) \to Q_\infty(x, y)$ as $n \to \infty$ for each $x, y \in \mathcal{S}$;
    \item[ii)] $Q_\infty$ is irreducible;
    \item[iii)] $\pi_n(x) \to \pi_\infty(x)$ as $n \to \infty$ for each $x \in \mathcal{S}$,
\end{enumerate}
then for each $x \in \mathcal{S}$,
\[
\sup_{t \geq 0} \| P_{n,x}(X(t) \in \cdot) - P_{\infty, x}(X(t) \in \cdot) \|_e \to 0
\]
$\text{as } n \to \infty$. 

\textbf{Proof.} We note that since $Q_\infty$ is non-explosive, irreducible, and positive recurrent (since it possesses a stationary distribution), 
\[
P_{\infty, x}(X(t) = y) \to \pi_\infty(y)
\]
as $t \to \infty$ for each $x\in \mathcal{S}$. Put $\lambda_n(x) = -Q_n(x,x)$ and let
\[
R_n(x, y) = 
\begin{cases} 
\frac{Q_n(x, y)}{\lambda_n(x)}, & x \neq y, \\ 
0, & x = y,
\end{cases}
\]
so that $\lambda_n(x)$ is the jump rate out of state $x \in \mathcal{S}$ and $R_n = (R_n(x, y) : x, y \in \mathcal{S})$ is the transition matrix of the embedded discrete-time Markov chain under $Q_n$. Given this,
\begin{equation}
    R_n(x, y) \to R_\infty(x, y) \label{eq:5.1}
\end{equation}
as $n \to \infty$ and $\lambda_n(x) \to \lambda_\infty(x)<\infty$  for $x \in \mathcal{S}$. Consequently, if $(Y_j:j\ge 0)$ is the embedded discrete-time Markov chain, the proof of Theorem 1 establishes that \eqref{eq:5.1} implies that
\begin{equation}
    P_{n,x}(Y_0 = y_0, \cdots, Y_m = y_m) \to P_{\infty,x}(Y_0 = y_0, \cdots, Y_m = y_m) \label{eq:5.2}
\end{equation}
as $n\to\infty$, for each $x \in \mathcal{S}$ and $n \geq 1$. Furthermore, if $T_i$ is the time of the $i$-th jump of $X$,
\[
T_i = \tau_0 +\cdots + \tau_{i-1},
\]
where $\tau_i$ is the time spent in the $i$-th state visited. Under $Q_n$, the $\tau_i$'s can be represented as
\[
\tau_i = \frac{1}{\lambda_n(Y_i)}\chi_i,
\]
where the $\chi_i$'s are iid exponential random variables with mean $1$. Since $\lambda_n(x) \to \lambda_\infty(x)$ as $n \to \infty$, evidently
\begin{equation}
    P_{n,x}((\tau_0, \tau_1, \dots, \tau_m) \in \cdot \mid Y_j, 0 \leq j \leq m) 
\Rightarrow 
P_{\infty,x}((\tau_0, \tau_1, \dots, \tau_m) \in \cdot \mid Y_j, 0 \leq j \leq m).\label{eq:5.3}
\end{equation}

Since
\[
X(t) = \sum_{m=0}^\infty Y_m \mathbb{I}(T_{m-1} \leq t < T_{m+1}),
\]
it follows from \eqref{eq:5.1}, \eqref{eq:5.2}, and \eqref{eq:5.3} that
\[
P_{n,x}(X(1) \in \cdot) \Rightarrow P_{\infty,x}(X(1) \in \cdot)
\]
as $n \to \infty$. We can then apply Corollary 1 to the sequence of transition matrices $(P_{n,x}(X(1) = y) : x, y \in \mathcal{S})$ to conclude that
\begin{equation}
    \sup_{m \geq 0} \| P_{n,x}(X(m) \in \cdot) - P_{\infty,x}(X(m) \in \cdot) \|_e \to 0
\quad \text{as } n \to \infty\label{eq:5.4}
\end{equation}
as $n\to\infty$. Since the $\mathcal{F}_e$ inclusion property ensures that
\[
\| P_{n,x}(X(t) \in \cdot) - \pi_n(\cdot) \|_e \to 0,
\]
is non-increasing in $t$, we obtain the theorem from \eqref{eq:5.4}.
$\Box$

The results of Section 3 are straightforward to extend to the jump process setting. The arguments of Theorem 2 and Proposition 3 go over without change in continuous time, upon recognizing that
\[
\pi_n(y) = \mathbb{E}_{n,x} \int_0^{\tau(x)} \mathbb{I}(X(t) = y) \, dt \quad \text{as } <\infty
\]
and
\[
\sum_x \pi_n(x) w(x) = \frac{\mathbb{E}_{n,x} \int_0^{\tau(x)} w(X(s)) \, ds}{\mathbb{E}_{n,x} \tau(x)},
\]
as $n\to\infty$, where $\tau(x)$ is the first time at which $X$ enters $x$ (i.e., $\tau(x) = \inf\{t \ge 0 : X(t) = x, X(t-) \neq x\}$) and $w(\cdot)$ is a non-negative function.  We obtain the following result.

\textbf{Proposition 5.} Suppose that $Q_\infty$ is irreducible and that $Q_n$ has a unique stationary distribution $\pi_n$ for $1\le n\le \infty$. If
\begin{enumerate}
    \item $Q_n(x, y) \to Q_\infty(x, y)$,
    \item $\pi_n(x) \to \pi_\infty(x)$,
    \item $\sum_{z \in \mathcal{S}} \pi_n(z) w(z) \to \sum_{z \in \mathcal{S}} \pi_\infty(z) w(z)$,
\end{enumerate}
as $n \to \infty$ for each $x, y \in \mathcal{S}$, then
\begin{align*}
 &   \mathbb{E}_{n,x} \int_0^T \exp\left(-\int_0^s \alpha(X(u)) \, du\right) r(X(s)) \, ds\\
 & \quad \to 
\mathbb{E}_{\infty,x} \int_0^T \exp\left(-\int_0^s \alpha(X(u)) \, du\right) r(X(s)) \, ds,
\end{align*}
for each non-negative $r: \mathcal{S} \to \mathbb{R}_+$ and $\alpha: \mathcal{S} \to \mathbb{R}_+$, for which $r \in \mathcal{F}_w$.

We conclude that stationary distribution convergence implies the strong interchange property, and also convergence of first transition expectations, in the Markov jump process setting.





\acks 
\noindent We thank the referees for their insightful comments, which led us to a sharper understanding of the key role of stationary distribution convergence in the current setting.



%
%
%
%

\bibliographystyle{APT}
\bibliography{ref}

\end{document}